\documentclass[11pt,aps,prb, preprint]{article}
\usepackage{amsmath}
\usepackage{amssymb}

\begin{document}
\title{Remarks on the Method of Modulus of Continuity and the Modified Dissipative Porous Media Equation}
\author{Kazuo Yamazaki}  
\date{}
\maketitle

\begin{abstract}
We employ Besov space techniques and the method of modulus of continuity to obtain the global well-posedness of the modified Porous Media Equation.\vspace{5mm}

\textbf{Keywords: Porous Media Equation, Quasi-geostrophic Equation, Criticality, Besov Space}
\end{abstract}
\footnote{2000MSC : 35B65, 35Q35, 35Q86}
\footnote{Department of Mathematics, Oklahoma State University, 401 Mathematical Sciences, Stillwater, OK 74078, USA, kyamazaki@math.okstate.edu}

\section{Introduction}
Porous Media Equation (PM) in $\mathbb{R}^{3}, t > 0$ is defined as follows:

\begin{equation}
\begin{cases}
\frac{\partial \theta}{\partial t} + u \cdot \nabla \theta + \nu \Lambda^{\alpha}\theta = 0; \hspace{5mm} \theta(x, 0) = \theta_{0}(x)\\
u = - \kappa(\nabla p + g\gamma\theta); \hspace{12mm} div u = 0
\end{cases}
\end{equation}

where $\theta$ represents liquid temperature, $\nu$ $>$ 0 the dissipative coefficient, $\kappa$ the matrix medium permeability divided by viscosity in different directions respectively, g the acceleration due to gravity, vector $\gamma$ the last canonical vector e$_{3}$ and $\Lambda = (-\triangle)^{1/2}$. Moreover, p is the liquid pressure and u represents the liquid discharge by Darcy's law. For simplicity we set $\kappa$ = g = 1.

In this paper we study the Modified Porous Media Equation (MPM) defined as follows:

\begin{equation}
\begin{cases}
\frac{\partial \theta}{\partial t} + u \cdot \nabla \theta + \nu \Lambda^{\alpha}\theta = 0; \hspace{5mm} \theta(x, 0) = \theta_{0}(x)\\
u = \Lambda^{\alpha - 1}\{-(\nabla p + \gamma\theta)\}; \hspace{5mm} div u = 0
\end{cases}
\end{equation}

in $\mathbb{R}^{3}$ with $\alpha \in (0, 1)$. Our main result is below:\\
\\
\textbf{Theorem 1.1}\\
\textit{Let $\nu >$ 0, 0 $< \alpha <$ 1 and $\theta_{0}$(x) $\in$ H$^{m}$, m $\in$ $\mathbb{Z}^{+}$, m $>$ 5/2. Then there exists a unique global solution $\theta$ to the MPM (2) such that}

\begin{equation*}
\theta \in \mathcal{C}(\mathbb{R}^{+}; H^{m}) \bigcap L^{2}_{loc}(\mathbb{R}^{+}; H^{m + \alpha /2}).
\end{equation*}

\textit{Moreover, for all $\gamma \in$ $\mathbb{R}^{+}$, we have t$^{\gamma} \theta \in$ L$^{\infty}(\mathbb{R}^{+}$, H$^{m + \gamma\alpha}$)}\\ 
\\
We note that an analogous version for Quasi-geostrophic Equation (QG), which we define below (3), was done by [\textbf{17}]. Moreover, while we will employ the method of Modulus of Continuity (MOC), initiated by [\textbf{12}] in a periodic setting, the smoothing effects stated above, i.e. the spatial decay of the solution, allows us to circumvent the difficulty in a non-periodic setting. In this regard, we cite [\textbf{1}] and [\textbf{9}] in which the authors proved the global well-posedness of QG with initial data belonging to the critical space $\dot{B}_{\infty, 1}^{0}$ and H$^{1}$ using the same technique. 
 
A similar result to the Theorem 1.1 showing global regularity of MPM (2) is also possible through the method introduced in [\textbf{5}] following the work in [\textbf{4}], [\textbf{6}] and [\textbf{7}]. A similar method following the work in [\textbf{13}] is also possible.

We stress that at first sight, modifying PM (1) by having $\Lambda^{\alpha - 1}$ act on the u term and finding its MOC based on the previous work on PM (1) in [\textbf{22}] seems somewhat difficult. As we will see, the u term of PM (1) can be decomposed to a linear combination of an identity and a singular integral operator acting on $\theta$ which we will denote by $\mathcal{P}(\theta)$. The problematic term is the Riesz potential, namely $\Lambda^{\alpha - 1}\theta$. We obviate from this issue by making a simple observation; see Proposition 3.3.

The outline of the rest of  the paper is as follows:
\begin{enumerate}
\item Introduction
\item Local Results
\item Global Results
\item Appendix A: Besov Space, Mollifers
\item Appendix B: Proofs of Local Results and More
\end{enumerate}

Let us introduce some MOC of relevance. By definition, a MOC is a continuous, increasing and concave function $\omega$: [0, $\infty$) $\to$ [0, $\infty$) with $\omega$(0) = 0. We say some function $\theta$: $\mathbb{R}^{n}$ $\to\mathbb{R}^{m}$ has MOC $\omega$ if $\lvert$$\theta$(x) - $\theta$(y)$\rvert$ $\leq$ $\omega$($\lvert$x - y$\rvert$) holds for all x, y $\in \mathbb{R}^{n}$.

The idea of MOC has caught much attention since the paper [\textbf{12}], in which the authors proved the global regularity of the solution to the 2-D critical QG defined as follows:

\begin{equation}
\partial_{t}\theta + (u \cdot \nabla) \theta = - \kappa\Lambda^{\alpha}\theta
\end{equation}

where u = (u$_{1}$, u$_{2}$) = (-$\mathcal{R}$$_{2}$$\theta$, $\mathcal{R}$$_{1}$$\theta$), $\mathcal{R}$$_{i}$ is Riesz Transform in $\mathbb{R}^{2}$, i = 1, 2, and $\kappa$ diffusivity constant. The variable u represents velocity and $\theta$ potential temperature. In particular, we have the following result from [\textbf{12}]:

\textbf{Proposition 1.2}\\
\textit{If the function $\theta: \mathbb{R}^{2}$ $\to \mathbb{R}$ has MOC $\omega$, then u of (3) has MOC as follows:}

\begin{equation*}
\Omega_{1} (\xi) = A (\int_{0}^{\xi}\frac{\omega(\eta)}{\eta}d\eta + \xi \int_{\xi}^{\infty}\frac{\omega(\eta)}{\eta^{2}}d\eta)
\end{equation*}

Their initiative motivated others to follow. Consider a pseudo-differential operator, or a modified Riesz Transform, $\tilde{\mathcal{R}}_{\alpha, j}$ defined as follows:

\begin{equation*}
\tilde{\mathcal{R}}_{\alpha, j} f(x) = \lvert D \rvert^{\alpha - 1}\mathcal{R}_{j}f(x) = c_{\alpha, n} \int_{\mathbb{R}^{n}}\frac{y_{j}}{\lvert y \rvert^{n + \alpha}}f(x - y)dy,
\end{equation*}

1 $\leq$ j $\leq$ n, for f $\in \mathcal{S}(\mathbb{R}^{n}$) where $\mathcal{R}$ is Riesz Transform, $\mathcal{S}$ Schwartz space, 0 $<$ $\alpha$ $<$ 1 and c$_{\alpha, n}$ the normalization constant. We have the following result due to [\textbf{17}]:

\textbf{Proposition 1.3}\\
\textit{If $\theta$, as defined in Proposition 1.2, has MOC $\omega$, then u = (-$\tilde{\mathcal{R}}_{\alpha, 2}\theta, \tilde{\mathcal{R}}_{\alpha, 1}\theta$) has the MOC of}

\begin{equation*}
\Omega_{2}(\xi) = A_{\alpha}(\int_{0}^{\xi}\frac{\omega(\eta)}{\eta^{\alpha}}d\eta + \xi\int_{\xi}^{\infty}\frac{\omega(\eta)}{\eta^{1+\alpha}}d\eta)
\end{equation*}

\textit{with some absolute constant A$_{\alpha} >$ 0 depending only on $\alpha$.}

Let us now derive a relation between u and $\theta$ of PM (1) following the method in [\textbf{3}]; for its generalization, see Lemma A.6 in the Appendix A. We have -curl (curl u) + $\nabla$ div u = $\triangle$u which is reduced to -curl(curl u) = $\triangle$u by divergence free property. Hence, we obtain

\begin{equation*}
\triangle u = (\frac{\partial^{2} \theta}{\partial x_{1} \partial x_{3}}, \frac{\partial^{2} \theta}{\partial x_{2} \partial x_{3}}, - \frac{\partial^{2}\theta}{\partial x_{1}^{2}} - \frac{\partial^{2}\theta}{\partial x_{2}^{2}}).
\end{equation*} 

Taking the inverse of the Laplacian,

\begin{equation*}
u = \frac{1}{4\pi} \int_{\mathbb{R}^{3}} \frac{1}{\lvert x-y \rvert}(\frac{\partial^{2} \theta}{\partial x_{1} \partial x_{3}}, \frac{\partial^{2} \theta}{\partial x_{2} \partial x_{3}}, - \frac{\partial^{2}\theta}{\partial x_{1}^{2}} - \frac{\partial^{2}\theta}{\partial x_{2}^{2}})dy
\end{equation*}

from which standard Integration By Parts (IBP) gives

\begin{equation*}
u(x, t) = - \frac{2}{3}(0, 0, \theta(x, t)) + \frac{1}{4\pi} PV \int_{\mathbb{R}^{3}}K(x - y)\theta(y, t)dy = C\theta + \mathcal{P}(\theta)
\end{equation*}

for x $\in \mathbb{R}^{3}$ where K(x) = ($\frac{3x_{1}x_{3}}{\lvert x \rvert^{5}}$, $\frac{3x_{2}x_{3}}{\lvert x \rvert^{5}}$, $\frac{2x_{3}^{2} - x_{1}^{2} - x_{2}^{2}}{\lvert x \rvert^{5}}$), considered as a kernel of double Riesz transform in [\textbf{22}]. For simplicity when no confusion arises, we let this constant C be one.

The scaling invariance of the PM is $\lambda$$^{\alpha - 1}$$\theta$ ($\lambda x, \lambda^{\alpha}t$) for $\lambda$ $>$ 0, same as the case of QG. This makes $\alpha$ = 1 the threshold of sub- and super-criticality. Recall that the method of MOC in the first place was introduced in order to prove the global regularity of the critical 2-D QG. Observe that if $\theta$(x, t) solves QG (3), then so does $\theta$($\lambda$x, $\lambda$t). The case when $\alpha < 1$ was studied in [\textbf{20}] and global regularity result was obtained under a certain initial condition. 

One may modify the QG (3) as below so that for any $\alpha \in$ (0, 1), its scaling invariance may be similar to that of the critical case (cf. [\textbf{4}]):

\begin{equation*}
\partial_{t}\theta + (u \cdot \nabla) \theta + \kappa \Lambda ^{\alpha}\theta = 0
\end{equation*}

with u = $\Lambda^{\alpha -1}(-\mathcal{R}_{2}\theta, \mathcal{R}_{1}\theta$). Observe that this PDE enjoys the rescaling of $\theta(x, t) \to \theta(\lambda x, \lambda^{\alpha}$t). For this reason, as we will see in (4) and (5), we will construct a MOC that is unbounded so that finding one MOC $\omega$ which is globally preserved in time implies that all the MOC $\omega_{\lambda}(\xi) = \omega(\lambda\xi)$ will also be globally preserved. It takes only a glance at MPM (2) to realize that it was defined in the same spirit.

\section{Local Results}
The purpose of this section is to state local results. We note that double Riesz Transform remains bounded in any space in which an ordinary Riesz Transform is bounded and that the results for the latter case has been obtained by [\textbf{17}]. The method of proof is similar to that described in [\textbf{14}] through regularizing (2) and relying on Picard's Theorem. Let us first set some notations; additional information can be found in Appendix A.

Denote by $\mathcal{S'}$ the space of tempered distributions, $\mathcal{S'}(\mathbb{R}^{n})/\mathcal{P}(\mathbb{R}^{n}$) the quotient space of tempered distributions modulo polynomials, $\mathcal{F}$f($\xi$) = $\hat{f}$($\xi$) the Fourier Transform, and $\lVert \cdot \rVert_{X}$ the norm of a Banach space X, e.g.

\begin{equation*}
H^{m} = \{f \in L^{2}(\mathbb{R}^{n}) : \lVert f \rVert_{H^{m}}^{2} = \sum_{0\leq\lvert\beta\rvert\leq m} \lVert D^{\beta}f \rVert_{L^{2}}^{2} < \infty\}.
\end{equation*}

Next we take the usual dyadic unity partition of Littlewood-Paley decomposition. Let us denote two nonnegative radial functions $\chi, \phi \in C^{\infty}(\mathbb{R}^{n}$) supported in \{$\xi \in \mathbb{R}^{n} : \lvert\xi\rvert \leq \frac{4}{3}\}$ and $\{\xi \in$ $\mathbb{R}^{n}$ : $\frac{3}{4}$ $\leq$ $\lvert\xi\rvert\leq\frac{8}{3}$\} respectively such that

\begin{equation*}
\chi(\xi) + \sum_{j\geq 0}\phi(2^{-j}\xi) = 1, \forall \xi \in \mathbb{R}^{n}; \hspace{5mm}  \sum_{j \in \mathbb{Z}}\phi(2^{-j}\xi) = 1,  \forall \xi \neq 0
\end{equation*}

and $\chi(\xi) = \hat{\Phi}(\xi), \phi(\xi) = \hat{\Psi}(\xi)$. We define for all $f \in\mathcal{S'} (\mathbb{R}^{n}$) the nonhomogeneous Littlewood-Paley operators:

\begin{equation*}
\triangle_{-1}f := \Phi \ast f, \triangle_{j}f := \Psi_{2^{-j}} \ast  f,  \forall j \in \mathbb{Z}^{+} \cup \{0\}
\end{equation*}

where $\widehat{\Psi_{2^{-j}}}(\xi) = \phi(2^{-j}\xi$) and the homogeneous defined as

\begin{equation*}
\dot{\Delta}_{j}f = \Psi_{2^{-j}} \ast f,  \forall j \in \mathbb{Z}.
\end{equation*}

With these Littlewood-Paley operators, we define Besov spaces for p, q $\in$ [1, $\infty$], s $\in \mathbb{R}$, the nonhomogeneous and homogeneous respectively:

\begin{equation*}
B^{s}_{p, r} := \{f \in \mathcal{S'}(\mathbb{R}^{n}): \lVert f \rVert_{B_{p, r}^{s}} := \lVert\triangle_{-1} f\rVert_{L^{p}} + (\sum_{j\geq 0}2^{jsr}\lVert\Delta_{j} f \rVert_{L^{p}}^{r})^{1/r} < \infty\}
\end{equation*}

\begin{equation*}
\dot{B}_{p,r}^{s} := \{f \in \mathcal{S'}(\mathbb{R}^{n})/\mathcal{P}(\mathbb{R}^{n}): \lVert f \rVert_{\dot{B}_{p,r}^{s}} := (\sum_{j \in \mathbb{Z}} 2^{jsr}\lVert \dot{\Delta}_{j} f \rVert_{L^{p}}^{r})^{1/r} < \infty\}.
\end{equation*}

Now let us state the main results of this section:

\textbf{Proposition 2.1}\\
\textit{Let $\nu >$ 0, 0 $< \alpha <$ 1 and $\theta_{0} \in$ H$^{m}$, m $\in \mathbb{Z}^{+}, m >$ 5/2. Then there exists a unique solution $\theta \in$ C([0, T], H$^{m}$) $\cap$ L$^{2}$([0, T], H$^{m+\frac{\alpha}{2}}$) to the MPM (2) where T $= T(\alpha, \lVert\theta_{0}\rVert_{m}) > 0$. Moreover, we have t$^{\gamma}\theta\in$ L$_{T}^{\infty}$H$^{m+\gamma\alpha}$ for all $\gamma \geq$ 1.}

\textbf{Proposition 2.2}\\
\textit{Let T$^{*}$ be the maximal local existence time of $\theta$ in C([0, T$^{*}$), H$^{m}$) $\cap$ L$^{2}$([0, T$^{*}$), H$^{m+\alpha/2}$). If T$^{*} < \infty$, then we have $\int_{0}^{T^{*}} \lVert\nabla\theta(t, \cdot)\rVert_{L^{\infty}}$dt = $\infty$.}

The proofs are found in the Appendix B.

\section{Global Results}
We extend our results to the global case. Below we use $\xi$ = $\lvert$x - y$\rvert$ interchangeably and omit the subscript $_{\alpha}$ but use instead i = 1, 2, 3, ... to indicate different constants; when no confusion arises. By the Blow-up Criterion Proposition 2.2, we only have to show that $\int_{0}^{T}$$\lVert\nabla\theta\rVert_{L^{\infty}}$ $<$ $\infty$ for $\theta$ a local solution of MPM (2) up to time T $>$ 0. We utilize the following two observations made in [\textbf{12}]:

\vspace{3mm}

\textbf{Proposition 3.1}\\
\textit{If $\omega$ is a MOC for $\theta(x, t): \mathbb{R}^{n}$ $\to \mathbb{R}$ for all t $>$ 0, then $\lvert\nabla\theta\rvert(x) \leq \omega'(0)$ for all x $\in \mathbb{R}^{n}$}

For this reason, we shall construct a MOC $\omega$ such that $\omega'(0) < \infty$. Next,

\textbf{Proposition 3.2}\\
\textit{Assume $\theta$ has a strict MOC satisfying}

\begin{equation*}
\omega'' (0 +) = - \infty
\end{equation*}

\textit{for all $t < T$; i.e. for all $x, y \in \mathbb{R}^{n}, \lvert\theta(x, t) - \theta(y, t)\rvert < \omega(\lvert x-y\rvert)$, but not for t $>$ T. Then, there exists x, y $\in \mathbb{R}^{n}$, x $\neq$ y such that $\theta$(x, T) - $\theta$(y, T) = $\omega$($\lvert$x - y$\rvert$).}

Consequently, the only scenario in which a MOC $\omega$ is lost is if there exists a moment $ T > 0$ such that $\theta$ has the MOC $\omega$ for all $ t \in [0, T]$ and two distinct points x and y such that $\theta$(x, T) - $\theta$(y, T) = $\omega$($\lvert$x - y$\rvert$). 

Below we rule out this possibility by showing that in such case, $\frac{\partial}{\partial t}$[$\theta$(x, t) - $\theta$(y, t)]$\lvert_{t = T}$ $<$ 0. For this purpose, let us write 

\begin{eqnarray*}
&&\frac{\partial}{\partial t}[\theta(x, T) - \theta(y, T)]\\
&=& - [( u \cdot \nabla \theta)(x, T) - (u \cdot \nabla \theta)(y, T)] - [(\Lambda^{\alpha} \theta)(x, T) - (\Lambda^{\alpha} \theta)(y, T)]
\end{eqnarray*}

We call the first bracket Convection term and the second Dissipation term. Our agenda now is to first estimate the Convection and Dissipation terms, to be specific find upper bounds that depend on $\omega$. Then we will construct the MOC $\omega$ explicitly that assures us that the sum of the two terms is negative to reach the desired result.

\textbf{Estimates on the Convection and Dissipation Terms}\\
We propose the following estimate for our Convection Term:

\textbf{Proposition 3.3}\\
\textit{If $\omega$ is a MOC for $\theta$(x,t): $\mathbb{R}^{3} \to \mathbb{R}$ , a local solution to MPM (2) for all $t \leq T$, then}

\begin{equation*}
(u \cdot \nabla \theta)(x, T) - (u \cdot \nabla \theta)(y, T) \leq C \Omega(\xi)\omega'(\xi)
\end{equation*}

\textit{where $\Omega(\xi$) is that of Proposition 1.3.}

The proof relies on the following observation due to [\textbf{15}] and [\textbf{19}]. From the expression of 

\begin{equation*}
-\triangle u = (-\partial_{x_{1}}\partial_{x_{3}}\theta, -\partial_{x_{2}}\partial_{x_{3}}\theta, \partial_{x_{1}}^{2}\theta + \partial_{x_{2}}^{2}\theta)
\end{equation*}

in the case of PM (1) derived in Section 1, the Fourier multilplier of such operator is clear; each component is a linear combination of terms like $\frac{\xi_{i}\xi_{j}}{\lvert\xi\rvert^{2}}, i, j = 1, 2, 3$ which belongs to $C^{\infty}(\mathbb{R}^{3}\setminus\{0\})$ and homogeneous of degree zero. Hence, it is clear that for the MPM (2) we can express u as

\begin{equation*}
\sum_{i, j}c_{ij}\Lambda^{\alpha-1}(-\triangle)^{-1}\partial_{i}\partial_{j}\theta = \sum_{i, j}\tilde{\mathcal{R}}_{\alpha, i, j}\theta
\end{equation*}

where $\tilde{\mathcal{R}}_{\alpha, i, j}$ may be considered as a modified double Riesz transform to which the result of Proposition 1.3 clearly applies. Therefore, u of MPM (2) has the same MOC of Proposition 1.3. 

\textbf{Remark}\\
We note that in the case of PM(1), making use of the observation above, we see that $\mathcal{P}(\theta)$ has $\Omega_{1}$ of Proposition 1.2 as its MOC; moreover, the additional term C$(\theta)$ can be absorbed into the first integral of $\Omega_{1}$. In [\textbf{22}], the authors considered the MOC of the Convection term of PM (1) separately; i.e. $\omega(\xi)$ for $\theta$ term and one different from $\Omega_{1}$ of Proposition 1.2 for the $\mathcal{P}(\theta)$. In our case, the same strategy would lead to having to deal with $\Lambda^{\alpha - 1}\theta$ term. 

Now we only need to compute

\begin{eqnarray*}
&& u \cdot \nabla \theta(x) - u \cdot \nabla \theta(y)\\
&=& \lim_{h \searrow 0}\frac{[\theta(x + hu(x)) - \theta(y + hu(y))] - [\theta(x) - \theta(y)]}{h}\\
&=& \lim_{h \searrow 0}\frac{[\theta(x + hu(x)) - \theta(y + hu(y))] - \omega(\xi)}{h}\\
&\leq& \lim_{h \searrow 0}\frac{\omega(\xi + h \Omega(\xi)) - \omega(\xi)}{h} = \Omega(\xi)\omega'(\xi)
\end{eqnarray*}

For the estimate on Dissipation term, we borrow below from [\textbf{20}], note the result is general in dimension:

\begin{equation*}
C_{2}[\int_{0}^{\xi/2}\frac{\omega(\xi + 2\eta) + \omega(\xi - 2\eta) - 2\omega(\xi)}{\eta^{1+\alpha}} d\eta + \int_{\xi/2}^{\infty}\frac{\omega(2\eta + \xi) - \omega(2\eta - \xi) - 2\omega(\xi)}{\eta^{1+\alpha}}d\eta]
\end{equation*}

\textbf{The Explicit Construction of the Modulus of Continuity}

We construct a relatively simple modulus of continuity. With $\alpha \in (0, 1), r \in (1, 1 + \alpha), 0 < \gamma < \delta (1-r\delta^{r-1})$, we define

\begin{equation}
\omega(\xi) = \xi - \xi^{r} \hspace{5mm} when \hspace{1mm} 0 \leq \xi \leq \delta
\end{equation}
\begin{equation}
\omega'(\xi) = \frac{\gamma}{\xi} \hspace{5mm} when \hspace{1mm} \delta < \xi
\end{equation}

The function $\omega$ is continuous and $\omega(0) = 0$. The first derivative of (4) is

\begin{equation}
\omega'(\xi) = 1 - r\xi^{r-1} \geq 0
\end{equation}

for $\delta$ sufficiently small and hence increasing . The unboundedness and $\omega'(0) < \infty$ can be readily checked. The second derivative of (4) is 

\begin{equation}
\omega''(\xi) = -r(r-1)\xi^{r-2} < 0
\end{equation}

From (5) we also have $\omega''(\xi) = -\frac{\gamma}{\xi^{2}} < 0$. Moreover, notice

\begin{equation}
lim_{\xi \to 0^{+}} \omega''(\xi) = -\infty
\end{equation}

Finally, 

\begin{equation}
\omega'(\delta^{+}) = \frac{\gamma}{\delta} < 1 - r\delta^{r-1} = \omega'(\delta^{-})
\end{equation}

Therefore, the concavity is achieved. Now we consider two cases

\textbf{The case $0 \leq \xi \leq \delta$}

\vspace{1mm}

In this case we make use of

\begin{equation}
\frac{\omega(\eta)}{\eta} \leq \omega'(0) = 1
\end{equation}

with which we immediately obtain

\begin{equation}
\int_{0}^{\xi}\frac{\omega(\eta)}{\eta^{\alpha}}d\eta \leq \int_{0}^{\xi}\frac{\omega(\eta)}{\eta}d\eta \leq \xi
\end{equation} 

since $\eta \leq \eta^{\alpha}$ for $\eta \leq \delta < 1$ and

\begin{equation}
\int_{\xi}^{\delta}\frac{\omega(\eta)}{\eta^{1 + \alpha}}d\eta = \int_{\xi}^{\delta} \frac{\eta - \eta^{r}}{\eta^{1 + \alpha}}d\eta \leq \int_{\xi}^{\delta}\eta^{-\alpha}d\eta \leq \frac{\delta^{1-\alpha}}{1-\alpha}
\end{equation}

Moreover,

\begin{eqnarray}
&&\int_{\delta}^{\infty}\frac{\omega(\eta)}{\eta^{1 + \alpha}}d\eta = \frac{\omega(\delta)\delta^{-\alpha}}{\alpha} + \frac{1}{\alpha}\int_{\delta}^{\infty}\omega'(\eta)\eta^{-\alpha}d\eta\\
&\leq& \frac{\delta^{1-\alpha}}{\alpha} + \frac{\gamma}{\alpha}\int_{\delta}^{\infty}\eta^{-\alpha - 1}d\eta =\frac{\delta^{1-\alpha}}{\alpha} + \frac{\gamma\delta^{-\alpha}}{\alpha^{2}} \leq (\frac{1}{\alpha} + \frac{1}{\alpha^{2}})\delta^{1-\alpha}\nonumber
\end{eqnarray}

where the first equality is by IBP and the last inequality used that $\gamma < \delta(1-r\delta^{r-1}) < \delta$. Since $\omega'(\xi) \leq \omega'(0) = 1$, the contribution from the positive side is limited to 

\begin{equation}
C_{1}[\xi + \xi[\frac{1}{1-\alpha} + \frac{1}{\alpha} + \frac{1}{\alpha^{2}}]\delta^{1-\alpha}]
\end{equation}

The work from [\textbf{20}] shows that the first integrand of the dissipation term gives 

\begin{equation}
\int_{0}^{\xi/2}\frac{\omega(\xi + 2\eta) + \omega(\xi - 2\eta) - 2\omega(\xi)}{\eta^{1 + \alpha}}d\eta \leq -C_{2}\xi^{r-\alpha}
\end{equation}

Therefore, adding (14) and (15) deduces the inequality of

\begin{equation}
\xi[C_{1} + C_{1}[\frac{1}{1-\alpha} + \frac{1}{\alpha} + \frac{1}{\alpha^{2}}]\delta^{1-\alpha} - C_{2}\xi^{r - \alpha - 1}] < 0
\end{equation}

to be achieved as $\delta \to 0$ forcing $\xi \to 0$; note $r - \alpha - 1 < 0$. 

\vspace{5mm}

\textbf{The Case $\delta < \xi$}

In this case we only have

\begin{equation}
\frac{\omega(\eta)}{\eta} \leq 1
\end{equation}

for $\eta \in [0, \delta]$. Using $\omega(\eta) \leq \omega(\xi)$ for all $\delta \leq \eta \leq \xi$, we obtain

\begin{eqnarray}
&&\int_{0}^{\xi}\frac{\omega(\eta)}{\eta^{\alpha}}d\eta = \int_{0}^{\delta}\frac{\omega(\eta)}{\eta^{\alpha}}d\eta + \int_{\delta}^{\xi}\frac{\omega(\eta)}{\eta^{\alpha}}d\eta\\
&\leq& \int_{0}^{\delta}\frac{\omega(\eta)}{\eta} + \omega(\xi)\int_{\delta}^{\xi}\frac{1}{\eta^{\alpha}}d\eta \leq \delta + \omega(\xi)\frac{\xi^{1-\alpha}}{1-\alpha}\nonumber
\end{eqnarray}

Observing that

\begin{equation}
\omega(\xi) \geq \omega(\delta) = \delta - \delta^{r} \geq \frac{\delta}{2}
\end{equation}

if $\delta$ is small enough, we deduce from (18)

\begin{equation}
\int_{0}^{\xi}\frac{\omega(\eta)}{\eta^{\alpha}}d\eta \leq \omega(\xi)(2 + \frac{\xi^{1-\alpha}}{1-\alpha})
\end{equation}

Next,

\begin{eqnarray}
&&\int_{\xi}^{\infty}\frac{\omega(\eta)}{\eta^{1 + \alpha}}d\eta =\frac{\omega(\xi)\xi^{-\alpha}}{\alpha} + \frac{1}{\alpha}\int_{\xi}^{\infty}\omega'(\eta)\eta^{-\alpha}d\eta\\
&=& \frac{\omega(\xi)\xi^{-\alpha}}{\alpha} + \frac{\gamma}{\alpha} \int_{\xi}^{\infty}\eta^{-1-\alpha}d\eta = \frac{\omega(\xi)\xi^{-\alpha}}{\alpha} + \frac{\gamma\xi^{-\alpha}}{\alpha^{2}} \leq \omega(\xi)\xi^{-\alpha}(\frac{\alpha + 1}{\alpha^{2}})\nonumber
\end{eqnarray}

if we take $\gamma$ small enough that 

\begin{equation}
2ln(2)\gamma < \frac{\delta}{2} \leq \omega(\xi) 
\end{equation}

following (19). Therefore, the contribution from the positive side is limited to

\begin{eqnarray}
&&C_{1}\Omega(\xi)\omega'(\xi) \leq C_{1}[\omega(\xi)(2 + \frac{\xi^{1-\alpha}}{1-\alpha}) + \omega(\xi)\xi^{1-\alpha}(\frac{\alpha + 1}{\alpha^{2}})]\omega'(\xi)\\
&=& C_{1}\omega(\xi)\xi^{-\alpha}[\frac{2(1-\alpha)\gamma\xi^{\alpha - 1} + \gamma}{(1-\alpha)} + \gamma(\frac{\alpha + 1}{\alpha^{2}})]\nonumber\\
&\leq& C_{1}\omega(\xi)\xi^{-\alpha}[\frac{2(1-\alpha)\gamma^{\alpha} + \gamma}{(1-\alpha)} + \gamma(\frac{\alpha + 1}{\alpha^{2}})]\nonumber
\end{eqnarray}

where we made use of $\alpha \in (0,1)$ and $\gamma < \delta < \xi$ in this case. On the other hand, we have the following estimate from [\textbf{20}]:

\begin{equation}
\int_{\xi/2}^{\infty}\frac{\omega(2\eta + \xi) - \omega(2\eta - \xi) - 2\omega(\xi)}{\eta^{1 + \alpha}}d\eta \leq -C_{2}\omega(\xi)\xi^{-\alpha}
\end{equation}

which is due to 

\begin{eqnarray*}
\omega(2\eta + \xi) - \omega(2\eta - \xi) \leq \omega(2\xi) &=& \omega(\xi) + \gamma\int_{\xi}^{2\xi}\frac{1}{\eta}d\eta = \omega(\xi) + \gamma ln(2) < \frac{3}{2}\omega(\xi)
\end{eqnarray*}

by concavity and (22). In sum, we obtain from (23) and (24) 

\begin{equation}
\omega(\xi)\xi^{-\alpha}[C_{1}[\frac{2(1-\alpha)\gamma^{\alpha} + \gamma}{(1-\alpha)} + \gamma(\frac{\alpha + 1}{\alpha^{2}})] - C_{2}] < 0
\end{equation}

for $\gamma$ sufficiently small. Q.E.D.

\section{Appendix A: Besov Space and Mollifiers}
\textbf{Besov Space}

We introduce the two types of coupled space-time Besov spaces. They are L$^{\rho}$([0, T], B$^{s}_{p,r}$), abbreviated by L$^{\rho}_{T}$B$^{s}_{p,r}$, defined with a norm

\begin{equation*}
\lVert f \rVert_{L_{T}^{\rho}B_{p, r}^{s}} = \lVert \lVert 2^{js}\rVert \triangle_{j}f\rVert_{L^{p}}\rVert_{l^{r}}\rVert_{L^{\rho}[0, T]}
\end{equation*}

and $\tilde{L}$$^{\rho}$([0, T], B$^{s}_{p,r}$), abbreviated by $\tilde{L}$$^{\rho}_{T}$B$^{s}_{p,r}$, called the Chemin-Lerner's space-time space, the set of tempered distribution f satisfying the norm

\begin{equation*}
\lVert f \rVert_{\tilde{L}_{T}^{\rho}B_{p,r}^{s}} := \lVert 2^{js}\lVert \Delta_{j} f \rVert_{L^{\rho}_{T}L^{P}}\rVert_{l^{r}} < \infty
\end{equation*}

We list useful results below:

\textbf{Lemma A.1 Bernstein's Inequality}\\
\textit{Let f $\in L^{p}(\mathbb{R}^{n})$ with 1 $\leq p \leq q \leq \infty$ and 0 $< r <$ R. Then for all $k \in \mathbb{Z}^{+}\cup\{0\}$, and $\lambda > 0$ there exists a constant C$_{k} >$ 0 such that}

\begin{equation}
\sup_{\lvert \alpha \rvert = k} \lVert\partial^{\alpha}f\rVert_{L^{q}} \leq C\lambda^{k+n(1/p - 1/q)}\lVert f\rVert_{L^{p}} \hspace{5mm} if supp \hspace{1mm}\mathcal{F}f \subset \{\xi : \lvert\xi\rvert \leq \lambda r\},
\end{equation}

\begin{equation}
C_{k}^{-1}\lambda^{k}\lVert f \rVert_{L^{p}} \leq sup_{\lvert\alpha\rvert = k} \lVert\partial^{\alpha}f\rVert_{L^{p}} \leq C_{k}\lambda^{k}\lVert f \rVert_{L^{p}} \hspace{5mm} if supp\mathcal{F}f \subset \{\xi : \lambda r \leq \lvert\xi\rvert \leq \lambda R\}.
\end{equation}

\textit{and if we replace derivative $\partial^{\alpha}$ by the fractional derivative, the inequalities remain valid only with trivial modifications.}\\
\\
\textbf{Lemma A.2 Besov Embedding (cf. [\textbf{18}])}\\
\textit{Assume s $\in \mathbb{R}$ and p, q $\in$ [1, $\infty$].}

\textit{(a) If 1 $\leq q_{1} \leq q_{2} \leq \infty$, then $\dot{B}_{p, q_{1}}^{s}(\mathbb{R}^{n}) \subset \dot{B}_{p, q_{2}}^{s}(\mathbb{R}^{n}$).}

\textit{(b) If 1 $\leq p_{1} \leq p_{2} \leq \infty$ and s$_{1} = s_{2} + n(\frac{1}{p_{1}} - \frac{1}{p_{2}})$, then $\dot{B}_{p_{1}, q}^{s_{1}}(\mathbb{R}^{n}) \subset \dot{B}^{s_{2}}_{p_{2}, q}(\mathbb{R}^{n})$.}

\textbf{Lemma A.3 (cf. [\textbf{21}])}

\textit{(a) For f $\in \mathcal{S}'$ with supp$\mathcal{F}f \subset \{\xi : \lvert\xi\rvert \leq r\}$, there exists C = C(n) such that for 1 $\leq p \leq q \leq \infty$,}

\begin{equation*}
\lVert f \rVert_{q} \leq Cr^{n(\frac{1}{p} - \frac{1}{q})}\lVert f \rVert_{p}
\end{equation*}

\textit{(b) Analogously, if supp$\mathcal{F}f \subset \{\xi : \lvert\xi\rvert \simeq r\}$, then}

\begin{equation*}
\lVert f \rVert_{q} \simeq r^{n(\frac{1}{p} - \frac{1}{q})}\lVert f \rVert_{p}
\end{equation*}

\textit{(c) Denoting Riesz transform by $\mathcal{R}$, for s $>$ n/p, 1 $< p < \infty$, $1 \leq r \leq \infty$,}

\begin{equation*}
\lVert\mathcal{R}f\rVert_{B_{p,r}^{s}} \leq C \lVert f\rVert_{B_{p,r}^{s}},
\end{equation*}

\textit{(d) Analogously, for 1 $\leq p \leq \infty$ and $1 \leq r \leq \infty$, we have}

\begin{equation*}
\lVert \mathcal{R}f\rVert_{\dot{B}_{p,r}^{s}} \leq C \lVert f\rVert_{\dot{B}_{p,r}^{s}}
\end{equation*}

Next, we define the transport-diffusion equation, 0 $< \alpha < 1$

\begin{equation*}
(TD)_{\alpha} 
\begin{cases}
\partial_{t}\theta + u \cdot \nabla \theta + \nu \lvert D \rvert^{\alpha} \theta = f; \hspace{5mm} \theta \lvert_{t=0} = \theta_{0}\\
div u = 0\\  
\end{cases}
\end{equation*}

\textbf{Proposition A.4 (cf. [\textbf{16}], [\textbf{17}])}

\textit{Let -1 $<$ s $<$ 1, 1$\leq\rho_{1}\leq\rho \leq \infty$, p, r $\in [1, \infty], f \in L^{1}_{loc}(\mathbb{R}^{+}, \dot{B}_{p,r}^{s+ \frac{\alpha}{\rho_{1}}-\alpha}$) and u a divergence-free vector field in L$^{1}_{loc}(\mathbb{R}^{+}$; Lip($\mathbb{R}^{n}$)). Suppose $\theta$ is a C$^{\infty}$ solution of (TD)$_{\alpha}$. Then there exists C = C(n, s, $\alpha$) such that for all t $\in\mathbb{R}^{+}$,}

\begin{equation}
\nu^{\frac{1}{\rho}}\lVert\theta\rVert_{\tilde{L}_{t}^{\rho}\dot{B}_{p,r}^{s+\alpha/\rho}} \leq Ce^{CV(t)}(\lVert\theta_{0}\rVert_{\dot{B}_{p,r}^{s}} + \nu^{\frac{1}{\rho_{1}}-1}\lVert f \rVert_{\tilde{L}_{t}^{\rho_{1}}\dot{B}_{p,r}^{s+\frac{\alpha}{\rho_{1}}-\alpha}}),
\end{equation}

\textit{where V(t) := $\int_{0}^{t}$$\lVert$$\nabla$u($\tau$)$\rVert_{L^{\infty}}$ d$\tau$.} 

The proof of Proposition A.4 consists of using para-differential calculus and Lagrangian coordinate method combined with commuter estimates; we refer readers to [\textbf{8}], [\textbf{10}], [\textbf{11}] and [\textbf{16}]. 

We also have the following from [\textbf{2}]:

\textbf{Proposition A.5}

\textit{Let u be a C$^{\infty}$ divergence-free vector field and f a C$^{\infty}$ function. Assume that $\theta$ is a solution of (TD$_{\alpha}$). Then for p $\in$ [1, $\infty$], we have}

\begin{equation*}
\lVert\theta(t)\rVert_{L^{p}} \leq \lVert\theta_{0}\rVert_{L^{p}} + \int_{0}^{t}\lVert f(\tau)\rVert_{L^{p}}d{\tau}
\end{equation*}

We introduce a result, relevant to our estimate of the Convection term:

\textbf{Lemma A.6 (cf. [19])}\footnote{The author is grateful to Professor Liutang Xue for pointing out this fact.}

\textit{Let m $\in C^{\infty}(\mathbb{R}^{n}\backslash\{0\})$ be a homogeneous function of degree 0, and T$_{m}$ be the corresponding multiplier operator defined by $\widehat{(T_{m}f)} = m\hat{f}$, then there exists a $\in \mathcal{C}$ and $\Omega \in C^{\infty}(\mathcal{S}^{n-1})$ with zero average such that for any Schwartz function f,}

\begin{equation*}
T_{m}f = af + PV\frac{\Omega(x')}{\lvert x \rvert^{n}}\ast f.
\end{equation*}

\textbf{Mollifier}

Given an arbitrary radial function $\rho(\lvert x \rvert) \in$ C$_{0}^{\infty}(\mathbb{R}^{3}), \rho \geq 0, \int_{\mathbb{R}^{3}}\rho$ dx = 1, we define the mollifer operator $\mathcal{T}_{\epsilon}$: L$^{p}(\mathbb{R}^{3}$)$\to$C$^{\infty}(\mathbb{R}^{3}), 1\leq p\leq \infty, \epsilon >$ 0, by

\begin{equation*}
(\mathcal{T}_{\epsilon}f)(x) = \epsilon^{-3}\int_{\mathbb{R}^{3}} \rho(\frac{x-y}{\epsilon})f(y)dy \hspace{1mm} \forall \hspace{1mm} f \in L^{p}(\mathbb{R}^{3}).
\end{equation*}

\textbf{Lemma A.7}

\textit{For m$\in \mathbb{Z}^{+}\cup\{0\}$, s $\in \mathbb{R}, k \in \mathbb{R}^{+}$. Then}

\textit{(a) For all $f \in C_{0}, \mathcal{T}_{\epsilon}f\to$f uniformly on a compact set U in $\mathbb{R}^{3}$; $\lVert\mathcal{T}_{\epsilon}f\rVert_{L^{\infty}}\leq\lVert f \rVert_{L^{\infty}}$.}

\textit{(b) For all $f\in H^{m}(\mathbb{R}^{3})$, D$^{\beta}(\mathcal{T}_{\epsilon}f$)=$\mathcal{T}_{\epsilon}(D^{\beta}$f); $\forall f \in H^{s}(\mathbb{R}^{3}), \lvert D \rvert^{s} (\mathcal{T}_{\epsilon}$f) = $\mathcal{T}_{\epsilon}(\lvert D\rvert^{s}$f).}

\textit{(c) For all $f \in H^{s}(\mathbb{R}^{3}), lim_{\epsilon \searrow 0}\lVert\mathcal{T}_{\epsilon}f - f\rVert_{H^{s}}$ = 0 and $\lVert\mathcal{T}_{\epsilon}f - f\rVert_{H^{s-1}}\leq c\epsilon\lVert f\rVert_{H^{s}}$.}

\textit{(d) For all  $f \in H^{m}(\mathbb{R}^{3}), \lVert\mathcal{T}_{\epsilon}f\rVert_{H^{m+k}} \leq c\epsilon^{-k} \lVert f \rVert_{H^{m}}$ and $\lVert\mathcal{T}_{\epsilon}D^{k}f\rVert_{L^{\infty}} \leq c \epsilon^{-(1+k)}\lVert f\rVert_{L^{2}}$.}

\section{Appendix B: Proofs of Local Results and More}
In this section we sketch the proofs from Section 2 Local Results. We regularize MPM (2) and study the following approximate system (ODE):

\begin{equation}
\begin{cases}
\theta_{t}^{\epsilon} + \mathcal{T}_{\epsilon}((\mathcal{T}_{\epsilon}u^{\epsilon})\cdot\nabla(\mathcal{T}_{\epsilon}\theta^{\epsilon})) + \nu\mathcal{T}_{\epsilon}(\mathcal{T}_{\epsilon}\Lambda^{\alpha}\theta^{\epsilon}) = 0; \hspace{5mm} \theta^{\epsilon}\lvert_{t=0} = \theta_{0}(x)\\
u^{\epsilon} = \Lambda^{\alpha -1}(C(\theta^{\epsilon}) + P(\theta^{\epsilon}))
\end{cases}
\end{equation}

Naturally the following Proposition can be proven by Picard's Theorem; we refer readers to [\textbf{14}] for proof.

\textbf{Proposition B.1  Global Existence of Regularized Solutions}

\textit{Let $\theta_{0}\in H^{m}, m\in \mathbb{Z}^{+}\cup\{0\},  m >$ 5/2. Then for all $\epsilon >$ 0, there exists  a unique global solution $\theta^{\epsilon} \in C^{1}([0, \infty), H^{m}$) to the regularized MPM (29).}

We may also assume the following with identical proof found in [\textbf{14}]:

\textbf{Proposition B.2}
\textit{The unique regularized solution $\theta^{\epsilon}\in C^{1}([0, \infty), H^{m}$) to (29) satisfies below:}

\begin{equation}
\frac{1}{2}\frac{d}{dt}\lVert\theta^{\epsilon}\rVert_{H^{m}}^{2} + \nu\lVert\mathcal{T}_{\epsilon}\Lambda^{\alpha/2}\theta^{\epsilon}\rVert_{H^{m}}^{2} \leq C_{m,\alpha}(\lVert\nabla\mathcal{T}_{\epsilon}\theta^{\epsilon}\rVert_{L^{\infty}} + \lVert\mathcal{T}_{\epsilon}\theta^{\epsilon}\rVert_{L^{3}})\lVert\theta^{\epsilon}\rVert_{H^{m}}^{2}
\end{equation}

\textit{and}

\begin{equation}
sup_{0 \leq t \leq T}\lVert\theta^{\epsilon}\rVert_{L^{2}} \leq \lVert\theta_{0}\rVert_{L^{2}}
\end{equation}

Below we show that there exists a subsequence convergent to a limit function $\theta$ that solves the MPM (2) up to some $T > 0$. The strategy is to first obtain the uniform bounds of H$^{m}$ norm in the interval [0, T] independent of $\epsilon$, and show that in [0, T], these approximate solutions are contracting in L$^{2}$ norm. By applying Interpolation Inequality, we will prove convergence as $\epsilon$ $\to$ 0 and pass the limit. Moreover, we outline the proof of the uniqueness and smoothing effects.

We first show below that ($\theta^{\epsilon}$) the family of solution is uniformly bounded in H$^{m}$. We have

\begin{equation*}
\frac{d}{dt}\lVert\theta^{\epsilon}\rVert_{H^{m}}\leq C_{m,\alpha}\lVert\mathcal{T}_{\epsilon}\nabla\theta^{\epsilon}\rVert_{L^{\infty}}\lVert\theta^{\epsilon}\rVert_{H^{m}} \leq C_{m, \alpha}\lVert\theta^{\epsilon}\rVert_{H^{m}}^{2}
\end{equation*}

by (30) and Sobolev Embedding as m $> \frac{5}{2}$. Thus, for all $\epsilon >$ 0, we have

\begin{equation}
sup_{0\leq t \leq T}\lVert\theta^{\epsilon}\rVert_{H^{m}}\leq\frac{\lVert\theta_{0}\rVert_{H^{m}}}{1-c_{m,\alpha}T\lVert\theta_{0}\rVert_{H^{m}}}
\end{equation}

which implies that for T $< \frac{1}{C_{m,\alpha}\lVert\theta_{0}\rVert_{H^{m}}}$, ($\theta^{\epsilon}$) is uniformly bounded in C([0, T], H$^{m}$). Next, by (30) and (32) after integrating in time [0, T] we obtain

\begin{equation}
\nu^{1/2}\lVert\Lambda^{\alpha/2}\theta^{\epsilon}\rVert_{L^{2}([0,T], H^{m})} \leq C(\lVert\theta_{0}\rVert_{H^{m}}, T)
\end{equation}

so that

\begin{equation*}
\nu^{1/2}\lVert\theta^{\epsilon}\rVert_{L^{2}([0,T], H^{m+\frac{\alpha}{2}})}\leq C(\lVert\theta_{0}\rVert_{H^{m}}, T), 
\end{equation*}

the desired uniform bound.

We now show that the solutions $\theta^{\epsilon}$ to regularized MPM (29) form a contraction in the low norm C([0,T], L$^{2}(\mathbb{R}^{3}$)); i.e. for all $\epsilon, \tilde{\epsilon},$ there exists C = C($\lVert\theta_{0}\rVert_{H^{m}}$, T) such that

\begin{equation*}
sup_{0 \leq t \leq T} \lVert\theta^{\epsilon} - \theta^{\tilde{\epsilon}}\rVert_{L^{2}}\leq C max \{\epsilon, \tilde{\epsilon}\}
\end{equation*}

We take

\begin{equation*}
\theta_{t}^{\epsilon}-\theta_{t}^{\tilde{\epsilon}} = -\nu(\mathcal{T}_{\epsilon}^{2}\Lambda^{\alpha}\theta^{\epsilon} - \mathcal{T}_{\tilde{\epsilon}}^{2}\Lambda^{\alpha}\theta^{\tilde{\epsilon}}) - [\mathcal{T}_{\epsilon}((\mathcal{T}_{\epsilon}u^{\epsilon})\cdot\nabla(\mathcal{T}_{\epsilon}\theta^{\epsilon}))-\mathcal{T}_{\tilde{\epsilon}}((\mathcal{T}_{\tilde{\epsilon}}u^{\tilde{\epsilon}})\cdot\nabla(\mathcal{T}_{\tilde{\epsilon}}\theta^{\tilde{\epsilon}}))]
\end{equation*}

and multiply by $\theta^{\epsilon} - \theta^{\tilde{\epsilon}}$ and integrate to get

\begin{eqnarray*}
&&(\theta_{t}^{\epsilon}-\theta_{t}^{\tilde{\epsilon}}, \theta^{\epsilon} - \theta^{\tilde{\epsilon}}) = -\nu(\mathcal{T}_{\epsilon}^{2}\Lambda^{\alpha}\theta^{\epsilon} - \mathcal{T}_{\tilde{\epsilon}}^{2}\Lambda^{\alpha}\theta^{\tilde{\epsilon}}, \theta^{\epsilon} - \theta^{\tilde{\epsilon}})\\
&-&([\mathcal{T}_{\epsilon}((\mathcal{T}_{\epsilon}u^{\epsilon})\cdot\nabla(\mathcal{T}_{\epsilon}\theta^{\epsilon})) - \mathcal{T}_{\tilde{\epsilon}}((\mathcal{T}_{\tilde{\epsilon}}u^{\tilde{\epsilon}})\cdot\nabla(\mathcal{T}_{\tilde{\epsilon}}\theta^{\tilde{\epsilon}}))], \theta^{\epsilon}-\theta^{\tilde{\epsilon}}) = I - II
\end{eqnarray*}

We bound I and II separately by standard method using the fact that Riesz potentials are bounded in L$^{p}$ space; for details, see [\textbf{17}]. Thus, we have

\begin{equation}
sup_{0\leq t \leq T}\lVert\theta^{\epsilon}-\theta^{\tilde{\epsilon}}\rVert_{L^{2}}\leq e^{c(M)T}(max\{\epsilon, \tilde{\epsilon}\} + \lVert\theta_{0}^{\epsilon} - \theta_{0}^{\tilde{\epsilon}}\rVert_{L^{2}}) \leq C(M, T) max\{\epsilon, \tilde{\epsilon}\}
\end{equation}

where M is an upper bound from (32). From this we deduce that $\{\theta^{\epsilon}\}$ is Cauchy in C([0, T], L$^{2}(\mathbb{R}^{3}$)) and hence converges to $\theta \in$ C([0, T], L$^{2}(\mathbb{R}^{3}$)). We apply the Interpolation Inequality to $\theta^{\epsilon} - \theta$, and using (32) and (34) we obtain

\begin{equation*}
sup_{0\leq t \leq T}\lVert\theta^{\epsilon}-\theta\rVert_{H^{s}} \leq C_{s}sup_{0 \leq t \leq T}(\lVert\theta^{\epsilon} - \theta\rVert_{L^{2}}^{1-\frac{s}{m}}\lVert\theta^{\epsilon}-\theta\rVert_{H^{m}}^{\frac{s}{m}})
\leq C(\lVert\theta_{0}\rVert_{H^{m}}, T, s)\epsilon^{1-\frac{s}{m}}
\end{equation*}

which gives $\theta \in C([0, T], H^{s}(\mathbb{R}^{3})), 0 \leq s < m$

Also, from $\theta_{t}^{\epsilon} = -\nu\mathcal{T}_{\epsilon}^{2}\Lambda^{\alpha}\theta^{\epsilon}$ - $\mathcal{T}_{\epsilon}((\mathcal{T}_{\epsilon}u^{\epsilon})\cdot\nabla(\mathcal{T}_{\epsilon}\theta^{\epsilon}$)), we see that $\theta_{t}^{\epsilon}$ converges to -$\nu\Lambda^{\alpha}\theta - u\cdot\nabla\theta$ in C([0, T], C($\mathbb{R}^{3}$)). As $\theta^{\epsilon} \to \theta$, the distribution limit of $\theta_{t}^{\epsilon}$ must be $\theta_{t}$; i.e. $\theta$ is a classical solution of MPM (2). From (32) and (33) we also have $\theta \in L^{\infty}([0, T], H^{m}(\mathbb{R}^{3}))\cap L^{2}([0,T], H^{m + \frac{\alpha}{2}}(\mathbb{R}^{3}$)).

Next, we show $\theta\in$ C([0, T], H$^{m}(\mathbb{R}^{3}$)). Firstly, we have 

\begin{eqnarray*}
&&\lVert\theta(t) - \theta(t')\rVert_{H^{m}}^{2} = C_{0}\lVert\theta(t) - \theta(t')\rVert_{B_{2,2}^{m}}^{2}\\
&\leq& C_{0}\sum_{-1 \leq j \leq J}2^{2jm}\lVert\triangle_{j}\theta(t) - \triangle_{j}\theta(t')\rVert_{L^{2}}^{2} + 2C_{0}\sum_{j > J}2^{2jm}\lVert\triangle_{j}\theta\rVert_{L_{T}^{\infty}L^{2}}^{2}
\end{eqnarray*}

and 

\begin{eqnarray*}
&&\lVert\nabla u\rVert_{L^{\infty}} + \lVert\Lambda^{\alpha}\theta\rVert_{L^{\infty}}\\
&\leq& C(\lVert\nabla\Lambda^{\alpha - 1}(\theta + \mathcal{P}(\theta))\rVert_{\dot{B}_{\infty, 1}^{0}} + \lVert\theta\rVert_{H^{m}})\\
&=& C(\sum_{j \leq -1}\lVert\dot{\triangle}_{j}\{\nabla\Lambda^{\alpha - 1}(\theta + \mathcal{P}(\theta))\}\rVert_{L^{\infty}} + \sum_{j \geq 0}\lVert\dot{\triangle}_{j}\{\nabla\Lambda^{\alpha - 1}(\theta + \mathcal{P}(\theta))\}\rVert_{L^{\infty}} + \lVert\theta\rVert_{H^{m}})\\
&\leq& C(\sum_{j \leq -1}2^{j \alpha}\lVert\dot{\triangle}_{j}\{\theta + \mathcal{P}(\theta)\}\rVert_{L^{\infty}} + \sum_{j \geq 0}2^{j(\alpha-1)}\lVert\dot{\triangle}_{j}\{\nabla(\theta + \mathcal{P}(\theta))\}\rVert_{L^{\infty}} + \lVert\theta\rVert_{H^{m}})\\
&\leq& C(\sum_{j \leq -1}2^{j(\alpha + \frac{3}{2})}\lVert\dot{\triangle}_{j}\theta\rVert_{L^{2}} + \sum_{j \geq 0}2^{j(\alpha -1)}\lVert\dot{\triangle}_{j}\nabla\theta\rVert_{L^{\infty}} + \lVert\theta\rVert_{H^{m}})\leq C\lVert\theta\rVert_{H^{m}}
\end{eqnarray*}

Next, by Besov embedding and Proposition A.4, we have

\begin{equation}
\lVert\theta\rVert_{\tilde{L}_{T}^{\infty}{\dot{B}}^{m}_{2,2}} \leq Ce^{cT\lVert\theta\rVert_{L_{T}^{\infty} H^{m}}}\lVert\theta_{0}\rVert_{H^{m}} < \infty
\end{equation}

Hence, we know there exists J = J(T, $\delta$) such that

\begin{equation*}
\sum_{j > J}2^{2jm}\lVert\triangle_{j}\theta\rVert_{L_{T}^{\infty}L^{2}}^{2} \leq \frac{\delta^{2}}{4C_{0}}
\end{equation*}

We apply Mean Value Theorem to get

\begin{eqnarray*}
&&\sum_{-1\leq j \leq J}2^{2jm}\lVert\triangle_{j}\theta(t) - \triangle_{j}\theta(t')\rVert_{L^{2}}^{2} \leq \lvert t - t'\rvert^{2} \sum_{-1 \leq j \leq J}2^{2jm}\lVert\triangle_{j}(\partial_{t}\theta)\rVert_{L_{T}^{\infty}L^{2}}^{2}\\
&\leq& C \lvert t - t'\rvert^{2}2^{2J}\lVert\partial_{t}\theta\rVert_{L_{T}^{\infty}H^{m-1}}^{2}
\end{eqnarray*}

On the last term, we have

\begin{eqnarray*}
&&\lVert\partial_{t}\theta\rVert_{H^{m-1}} \leq \nu \lVert\Lambda^{\alpha}\theta\rVert_{H^{m-1}} + \lVert u\cdot\nabla\theta\rVert_{H^{m-1}}\\
&\leq& \nu \lVert\theta\rVert_{H^{m-1+\alpha}} + \lVert u\theta\rVert_{H^{m}} \leq C(\lVert\theta\rVert_{H^{m}} + \lVert u \rVert_{H^{m}} \Vert\theta\rVert_{H^{m}})\\
&\leq& C(\lVert\theta\rVert_{H^{m}} + \lVert\theta\rVert_{H^{m}}^{2}) \leq C(\lVert\theta_{0}\rVert_{H^{m}}, T)
\end{eqnarray*}

i.e. $\partial_{t}\theta \in L^{\infty}([0,T], H^{m-1}$); hence the desired continuity.

The uniqueness is proven by standard way of using the difference of two different solutions, multiplication, integration and Gronwall's inequality (cf. [\textbf{17}]). For the smoothing effects, take t$^{\gamma}\theta, \gamma >$ 0 in (TD)$_{\alpha}$ below:

\begin{equation*}
\partial_{t}(t^{\gamma}\theta) + u\cdot\nabla(t^{\gamma}\theta) + \nu\Lambda^{\alpha}(t^{\gamma}\theta) = \gamma t^{\gamma-1}\theta;  (t^{\gamma}\theta)\rvert_{t=0} = 0
\end{equation*}

Assume T $\geq$ 1 without loss of generality. We show the following:

\begin{equation*}
\lVert t^{\gamma}\theta(t)\rVert_{L_{T}^{\infty}L^{2}} + \lVert t^{\gamma}\theta(t)\rVert_{\tilde{L}_{T}^{\infty}\dot{H}^{m+\gamma\alpha}}\leq C(T^{\gamma+1} + e^{c(\gamma+1)T\lVert\theta\rVert_{L_{T}^{\infty}H^{m}}})\lVert\theta_{0}\rVert_{H^{m}}
\end{equation*}

which implies 

\begin{equation*}
\lVert t^{\gamma}\theta(t)\rVert_{L_{T}^{\infty}H^{m+\gamma\alpha}}\leq C(T^{\gamma+1} + e^{C(\gamma+1)T\lVert\theta\rVert_{L_{T}^{\infty}H^{m}}})\lVert\theta_{0}\rVert_{H^{m}}
\end{equation*}

The proof is done through induction on $\gamma$ and interpolation to apply for all $\gamma \in \mathbb{R}^{+}$; the readers are referred to [\textbf{17}] for detail.

Finally, the blow up criterion is proven. In similar fashion to (30) we can obtain

\begin{equation*}
\frac{1}{2}\frac{d}{dt}\lVert\theta\rVert_{H^{m}}^{2} + \nu\lVert\Lambda^{\frac{\alpha}{2}}\theta\rVert_{H^{m}}^{2} \leq C_{m, \alpha}(\lVert\nabla\theta\rVert_{L^{\infty}} + \lVert\theta\rVert_{L^{3}})\lVert\theta\rVert_{H^{m}}^{2}
\end{equation*}

Gronwall's inequality shows that if the blow-up time T$^{\ast} < \infty$, then

\begin{equation*}
\int_{0}^{T^{\ast}}\lVert\nabla\theta\rVert_{L^{\infty}} dt = \infty
\end{equation*}

This completes the proofs of both Proposition 2.1 and 2.2.

\section{Acknowledgment}
The author expresses gratitude to Professor Gautam Iyer and Jiahong Wu for their teaching and Professor Miao Changxing, Lenya Ryzhik and Liutang Xue and the referee for their helpful comments and suggestions.

\end{document}